\documentclass[12pt,a4paper]{amsart}
\usepackage{amssymb}
\usepackage{mathrsfs}

\headsep=1cm \numberwithin{equation}{section}
\newtheorem{thm}{Theorem}[section]
\newtheorem{cor}[thm]{Corollary}
\newtheorem{lem}[thm]{Lemma}
\newtheorem{prop}[thm]{Proposition}
\theoremstyle{definition}
\newtheorem{defn}[thm]{Definition}
\theoremstyle{remark}
\newtheorem{rem}[thm]{Remark}
\numberwithin{equation}{section}
\newcommand{\nat}{\mathbb N}
\newcommand{\inte}{\mathbb Z}
\newcommand\Supp{\operatorname{Supp}}
\newcommand\Ass{\operatorname{Ass}}

\newcommand\Spec{\operatorname{Spec}}

\newcommand\Hom{\operatorname{Hom}}
\newcommand\Ext{\operatorname{Ext}}

\newcommand{\p}{\frak p }
\newcommand{\q}{\frak q }
\begin{document} \title [Faltings' local-global principle  for the finiteness of local cohomology]
{Faltings' local-global principle  for the finiteness of local cohomology modules over Noetherian rings}
\author[Ali Akbar Mehrvarz, Reza Naghipour and Monireh Sedghi]{Ali Akbar Mehrvarz$^*$, Reza Naghipour and Monireh Sedghi}
%\vspace*{0.5cm}Dedicated to Professor Hossein Zakeri}
\address{Department of Mathematics, Tabriz branch, Islamic Azad University, Tabriz, Iran}
\email{amehrvarz2013@gmail.com}
\address{Department of Mathematics, University of Tabriz, Tabriz, Iran.}
\email{naghipour@ipm.ir} \email {naghipour@tabrizu.ac.ir}

\address{Department of Mathematics, Azarbaijan Shahid Madani University, Tabriz, Iran. }%
\email{sedghi@azaruniv.ac.ir}  \email {m\_sedghi@tabrizu.ac.ir}%
\thanks{ 
%This research of the second author was  in part supported by a grant from IPM\\
$^*$Corresponding author: e-mail: {\it amehrvarz2013@gmail.com} (Ali Akbar Mehrvarz)}%
%\subjclass{}%

%\date{}%
%\dedicatory{}%
%\commby{}%
% ----------------------------------------------------------------
\begin{abstract}
Let $R$ denote a commutative Noetherian (not necessarily local) ring, $\frak a$ an ideal of $R$ and $M$ a finitely generated $R$-module.
The purpose of this paper is to  show that
$f^n_{\frak a}(M)={\rm inf} \{0\leq i\in\inte|\, \dim H^{i}_{\frak a}(M)/N \geq n  \text{for any finitely generated submodule}\,\, N \subseteq H^{i}_{\frak a}(M)\}$,
where $n$ is  a non-negative integer and the invariant  $f^n_{\frak a}(M):=\inf\{f_{\frak a R_{\frak p}}(M_{\frak p})\,\,|\,\,{\frak p}\in \Supp M/\frak a M\,\,{\rm and}\,\,\dim R/{\frak p}\geq n\}$
is  the $n$-th finiteness dimension of $M$ relative to $\frak a$. As a consequence, it follows that the set
$$ \Ass_R(\oplus _{i=0}^{f^n_{\frak a}(M)}H^{i}_{\frak a}(M))\cap \{\frak p\in \Spec R|\, \dim R/\frak p\geq n\}$$ is finite. This generalizes the main result of Quy \cite{Qu}, Brodmann-Lashgari \cite{BL} and Asadollahi-Naghipour \cite{AN}.
\end{abstract}
\maketitle
% ----------------------------------------------------------------
\section{Introduction}
Let $R$ be a commutative Noetherian ring, $\frak a$ an ideal of $R$, and
$M$ a finitely generated $R$-module. An important theorem
in local cohomology is Faltings${^,}$ Local-global Principle for the
finiteness of local cohomology modules (\cite[Satz 1]{F}) which
states that for a positive integer $r$ the $R_{\frak p}$-module
$H_{\frak a R_{\p}}^i(M_{\p})$ is finitely generated for all $i< r$ and for
all ${\p}\in \Spec R$ if and only if the $R$-module $H_{\frak a}^i(M)$ is
finitely generated for all $i< r$.

Another formulation of
Faltings${^,}$ Local-global Principle, particularly relevant for
this paper, is in terms of the finiteness dimension $f_{\frak a}(M)$ of $M$
relative to $I$, where
$$f_{\frak a}(M):=\inf\{i\in\mathbb{N}: H_{\frak a}^i(M)\text{ is not finitely generated}\},$$
with the usual convention that the infimum  of the empty set of
integers is interpreted as $\infty$.
 Bahmanpour et al., in  \cite{BNS}, introduced the notion of the $n$-th finiteness dimension
$f_{\frak a}^n(M)$ of $M$ relative to $I$ by
$$f_{\frak a}^n(M):=\inf\{f_{{\frak a}R_{\p}}(M_{\p}): {\p}\in\Supp M/{\frak a}M\text{ and} \dim R/{\p}\geq n\}.$$
More recently,  Asadollahi and Naghipour in \cite{AN},  introduced the class of
{\it in dimension $<n$ modules}, and they showed that,  if $(R,{\frak m})$ is a complete local ring, ${\frak a}$  an ideal of $R$ and $M$ a finitely generated $R$-module, then for any $n\in\mathbb{N}_0$,
$$f_{\frak a}^n(M)={\rm inf} \{0\leq i\in\inte|\, H^{i}_{\frak a}(M)\, \text{is not in dimension}\, < n\}. $$

In this paper, we eliminate the complete local hypothesis entirely  by proving the following:

\begin{prop}
Let $R$ be a  Noetherian ring, ${\frak a}$  an ideal of $R$  and $M$ a finitely generated $R$-module. Then for any $n\in\mathbb{N}_0$,
$$f_{\frak a}^n(M)={\rm inf} \{0\leq i\in\inte|\, H^{i}_{\frak a}(M) \text{is not in dimension}\, < n\}. $$
\end{prop}

The proof of Proposition 1.1  is given in Theorem 2.10. Pursuing this point of view further we establish the following consequence of Proposition 1.1.

\begin{thm}
Let $R$ be a Noetherian ring and ${\frak a}$ an ideal of $R$.   Let $n$ be a non-negative integer and $M$ a
finitely generated $R$-module. Then the $R$-modules  $\Ext_R^j(R/{\frak a}, H_{\frak a}^i(M))$ are in dimension $<n$ for all $i< f_{\frak a}^n(M)$ and all integers $j$. Moreover,  the $R$-modules  $\Ext_R^j(N, H_{\frak a}^{f_{\frak a}^n(M)}(M))$ are in dimension $<n$, for each finitely generated $R$-module $N$ with support in $V(\frak a)$ and $j=0, 1$.
\end{thm}

As a consequence of  Proposition 1.1  and Theorem 1.2, we derive the following, which is a generalization of the main result of Quy \cite[Theorem 3.2]{Qu}  and  Brodmann-Lashgari \cite[Theorem 2.2]{BL}.

\begin{cor}
Let $R$ be a Noetherian ring, ${\frak a}$  an ideal of $R$ and $M$ a finitely generated $R$-module. Then for any $n\in\mathbb{N}_0$,
the set
%\begin{center}
$$\Ass_R(\oplus _{i=0}^{f^n_{\frak a}(M)}H^{i}_{\frak a}(M))\cap \{\frak p\in \Spec R|\, \dim R/\frak p\geq n\},$$
%\end{center}
is finite.
\end{cor}
Another consequence of Proposition 1.1 shows that the first non-minimax (resp. non-skinny)  local cohomology modules
$H^{i}_{\frak a}(M)$ of a finitely module $M$ over a  Noetherian ring $R$ with respect to an ideal $\frak a$ is
$f_{\frak a}^n(M)={\rm inf} \{0\leq i\in\inte|\, H^{i}_{\frak a}(M)\, \text{is not in dimension}\, < 1\}$ (resp.
$f_{\frak a}^n(M)={\rm inf} \{0\leq i\in\inte|\, H^{i}_{\frak a}(M)\, \text{is not in dimension}\, < 2\}$).\\

Throughout this paper, $R$ will always be a commutative Noetherian ring with non-zero identity and $\frak a$ will be an ideal of $R$.
For an $R$-module $L$, the $i$-th local cohomology module of $L$ with support in $V(\frak a)$
is defined as:
$$H^i_{\frak a}(L) = \underset{n\geq1} {\varinjlim}\,\, {\rm Ext}^i_R(R/\frak a^n, L).$$

Local cohomology was defined and studied by Grothendieck. We refer the reader to \cite{BS} or \cite{Gr1} for more
details about local cohomology.

For a non-negative integer $n$,  an $R$-module $M$ is said to be {\it  in dimension $<n$} if there is a finitely generated submodule $N$ of $M$ such that $\dim M/N<n$.   By a {\it skinny or weakly
Laskerian} module, we mean an $R$-module $M$ such that the set $\Ass_R M/N$ is finite, for each submodule $N$ of $M$  (cf. \cite{Ro} or \cite{DM1}). Moreover, an $R$-module $M$ is said to be {\it minimax}, if there exists a finitely generated submodule $N$ of $M$, such that $M/N$ is Artinian. The class of minimax modules was introduced by H. Z\"{o}schinger \cite{Zo1} and he has given in \cite{Zo1, Zo2} many equivalent conditions for a module to be minimax.

\section{Main  Results}

In \cite{Qu} , P. H. Quy introduced the class of FSF modules and he has given some properties of this modules. The $R$-module $M$ is said to be a FSF module if there is a finitely generated submodule $N$ of $M$ such that the support of the quotient module $M/N$ is finite. When $R$ is a Noetherian ring, it is clear that, if $M$ is FSF, then $\dim \Supp M/N\leq 1$. This motivates the following definition.

%%%%%%%%%%%%%%%%%%%%%%%%%%%%%%%definition%%%%%%%%%%%%%%%%%%%%%%%%%%%%%%%%%%%%%%%%%%%%%%%%%%%%%%%%%%%%%%%%%%%%%%%%%%%%%%%%%%%%%%%%%%%%%%%%%%%%%
\begin{defn}
Let $n$ be a non-negative integer. An $R$-module $M$ is said to be in dimension $< n$, if there is a finitely generated submodule $N$ of $M$ such that
$\dim \Supp  M/ N <n$.
\end{defn}

\begin{rem}\label{rem}
Let $n$ be a non-negative integer and let $M$ be an $R$-module.
   \begin{enumerate}
     \item
     if $n=0$, then $M$ is in dimension $< n$ if and only if $M$ is Noetherian.
     \item
     If $M$ is minimax, then $M$ is in dimension $< 1$. In particular, if $M$ is Noetherian or Artinian, then $M$ is in dimension $< 1$.
     \item
     If $M$ is FSF, then $M$ is in dimension $< 2$.
     \item
     If $M$ is skinny, then $M$ is in dimension $< 2$, by \cite [Theorem 3.3]{Ba}.
     \item
     If $M$ is reflexive, then $M$ is in dimension $<1$.
     \item
     If $M$ is linearly compact, then $M$ is in dimension $<1$.
  \end{enumerate}
\end{rem}

\begin{defn}
If $T$ is an arbitrary subset of $\Spec R$ and $n\in\nat_0$, then we set $$(T)_{\geq n}:=\{ \p\in T | \dim R/\p\geq n\}.$$
\end{defn}

\begin{defn}
Let $R$ be a Noetherian ring,  $\frak a$  an ideal of $R$ and $M$ an $R$-module. For any  non-negative integer $n$, we define
$${\rm h}_{\frak a}^n(M):={\rm inf} \{0\leq i\in\inte: \,\, H^{i}_{\frak a}(M) \,\, \text{is not in dimension}< n\}. $$
\end{defn}

\begin{lem}\label{lem1}
Let $R$ be a Noetherian ring,  $n$  a non-negative integer and $L$ an $R$-module such that $\dim L\geq n.$ If there is a submodule $L'$ of $L$ such that $\dim L/L'\geq n$, then

\begin{center}
$\cap _{\frak p \in (\Ass_R  L)_{\geq n}} \frak p\subseteq \cap_{\frak p \in (\Ass_R  L/L')_{\geq n}} \frak p.$
\end{center}
%$$\{\frak p\in \Ass_RL: \dim R/\frak p\geq s\}\subseteq  \{\frak p\in \Ass_R L/L': \dim R/\frak p\geq s\}.$$
\end{lem}

\proof  Let $\p \in (\Ass_R L/L')_{\geq n}$. Then $\p\in\Supp L$, and so there exists $\q\in\Ass_R L$ such that $\q\subseteq \p$. Therefore $\q\in(\Ass_R L)_{\geq n}$, and thus

\begin{center}
$\cap_{\p \in (\Ass_R L)_{\geq n}}\frak p\subseteq \q\subseteq \p.$
\end{center}
This completes the proof.\qed \\
%%%%%%%%%%%%%%%%%%%%%%%%%%%%%%%%%%%%%%%%%%%%lemma%%%%%%%%%%%%%%%%%%%%%%%%%%%%%%%%%%%%%%%%%%%%%%%%%%%%%%%%%%%%%%%%%%%%%%%%%%%%%%%%%%%%%%%%%%%%%%%%%%%%%%
\begin{lem}\label{lem2}
Let $R$ be a Noetherian ring,  $n$  a non-negative integer and $L$  an $R$-module in dimension $<n$. Then the set $(\Ass_R L)_{\geq n}$ is finite.
\end{lem}

\proof
Since $L$ is in dimension $< n$, it follows from the definition that there is a finitely generated submodule $L'$ of $L$ such that $\dim \Supp L/L'<n$. Now, from the exact sequence $$0 \longrightarrow L' \longrightarrow L\longrightarrow L/L' \longrightarrow 0$$
we obtain
$$(\Ass_R L)_{\geq n}\subseteq (\Ass_R L')_{\geq n}\cup (\Ass_R L/L')_{\geq n}.$$
As $\dim\Supp L/L'<n$, it follows that $(\Ass_R L/L')_{\geq n}=\emptyset$. Thus  $$(\Ass_R L)_{\geq n}\subseteq(\Ass_R L')_{\geq n},$$ and so  the set $(\Ass_R L)_{\geq n}$ is finite.\qed \\

\begin{lem}
Let $R$ be a Noetherian ring, $s$  a non-negative integer, $M$  a finitely generated $R$-module and $\frak a$ an ideal of $R$. Let ${\bf x}=x_1,\hdots, x_s\in {\frak a}$ be an $M$-regular sequence. Then
$$f_{\frak a}(M)\leq s+f_{\frak a}(M/{\bf x}M).$$
\end{lem}

\proof  The assertion follows easily by induction on $s$.\qed \\

\begin{cor}
\label{1}
Let $R$ be a Noetherian ring, $s$  non-negative integer, $M$  a finitely generated $R$-module and ${\frak a}$ an ideal of $R$.  Let  ${\bf x}=x_1,\hdots, x_s\in {\frak a}$ be an $M$-regular sequence. Then for any $n\in \mathbb{N}_0$,
$$f_{\frak a}^n(M)\leq s+f_{\frak a}^n(M/{\bf x}M).$$
\end{cor}
\proof
The assertion follows from  the definition of $f_{\frak a}^n(M)$ and  induction on $s$.\qed \\

The following proposition which plays a key role in this paper will serve to shorten the proof of the main theorems.

\begin{prop}
\label{prop}
Let $R$ be a Noetherian ring,  $n$  a non-negative integer, $M$ a finitely generated $R$-module and ${\frak a}$ an ideal of $R$. Then for all $i<f_{\frak a}^n(M)$, the $R$-module $H_{\frak a}^i(M)$ is in dimension $< n$.
\end{prop}

\proof
We show that $H_{\frak a}^i(M)$ is in dimension $< n$  by induction on $i$.
The case $i=0$ is clear, because $H_{\frak a}^0(M)$ is finitely generated.
So suppose that $i>0$ and that the result has been proved for smaller values of $i$. By this inductive assumption,
$H_{\frak a}^j(M)$ is in dimension $< n$  for $j=0,1, \dots, i-1$, and it only remains for us to prove that
$H_{\frak a}^i(M)$ is in dimension $< n$. To this end, it follows from \cite[Corollary 2.1.7 and Lemma 2.1.2]{BS}  that
$H_{\frak a}^i(M)\cong H_{\frak a}^i(M/H^0_{\frak a}(M))$ and $M/H^0_{\frak a}(M)$ is an $\frak a$-torsion-free $R$-module. Hence
we can (and do) assume that $M$ is  an $\frak a$-torsion-free $R$-module.

We now use \cite[Lemma 2.1.1]{BS} to deduce that  ${\frak a}$ contains an element $x$ which is $M$-regular. Let $t\in\nat$.  Then the exact sequence
$$0 \longrightarrow M \stackrel{x^t} \longrightarrow  M \longrightarrow M/x^tM
\longrightarrow 0$$ induces a long exact sequence
$$\cdots \longrightarrow H^{i}_{{\frak a}}(M)\stackrel{x^t} \longrightarrow H^{i}_{{\frak a}}(M)
\longrightarrow H^{i}_{{\frak a}}(M/x^tM) \longrightarrow
H^{i+1}_{{\frak a}}(M)\stackrel{x^t} \longrightarrow H^{i+1}_{{\frak a}}(M)
\longrightarrow \cdots,$$
and so we obtain  the exact sequence

$$0 \longrightarrow H^{i-1}_{{\frak a}}(M)/x^tH^{i-1}_{{\frak a}}(M) \longrightarrow H^{i-1}_{{\frak a}}(M/x^tM)
\longrightarrow (0 :_{H^i_{{\frak a}}(M)} x^t) \longrightarrow 0.\,\,\,\,\,\,\,\,\,\,\,\,\, (\dag) $$

Since by Corollary \ref{1}, $$f_{\frak a}^n(M)\leq 1+f_{\frak a}^n(M/x^tM),$$
it follows from the inductive hypothesis that the $R$-module $H_{\frak a}^j(M/x^tM)$ is in dimension $< n$, for all $0\leq j\leq i-1$.

 Now, in order to show that $H_{\frak a}^i(M)$ is in dimension $< n$, suppose the contrary is true. Then, as  $H_{\frak a}^{i-1}(M/x^tM)$  is in dimension $<n$, it follows from the exact sequence  ($\dag$)  that the $R$-module $(0:_{H_{\frak a}^i(M)}x^t)$ is in dimension $<n$. Therefore, in view of Lemma \ref{lem2} the set $(\Ass_R(0:_{H_{\frak a}^i(M)}x^t))_{\geq n}$ is finite. Consequently  the set
$(\Ass_R H_{\frak a}^i(M))_{\geq n}$ is also finite.
Let $$(\Ass_R H_{\frak a}^i(M))_{\geq n}=\{\p_1,\hdots,\p_r\}.$$
Then, in view of the definition, the $R_{\p_j}$-module $(H_{\frak a}^i(M))_{\p_j}$ is finitely generated for all $1\leq j\leq r$. Thus for all $1\leq j\leq r$, there exists a finitely generated submodule $N_j$ of $H_{\frak a}^i(M)$ such that $(H_{\frak a}^i(M))_{\p_j}=(N_j)_{\p_j}$.
Set $K_{1}=N_1+\cdots +N_r$. Then $K_{1}$ is a finitely generated submodule of $H_{\frak a}^i(M)$ and we have
$$(\Ass_R H_{\frak a}^i(M)/K_{1})_{\geq n}\cap (\Ass_R H_{\frak a}^i(M))_{\geq n}=\emptyset.$$
Since $K_{1}$ is a finitely generated submodule of $H_{\frak a}^i(M)$, it follows that there exists a non-negative integer $l$ such that $K_{1}\subseteq (0:_{H_{\frak a}^i(M)}x^l)$, and so $$(K_{1}:_{H_{\frak a}^i(M)}x)\subseteq (0:_{H_{\frak a}^i(M)}x^{l+1}).$$

Now as the $R$-module $(0:_{H_{\frak a}^i(M)}x^{l+1})$ is in dimension $< n$, it follows that the $R$-module $(K_{1}:_{H_{\frak a}^i(M)}x)/K_{1}$ is in dimension $< n$. Therefore, using again the above method, we see that the set $(\Ass_R (H_{\frak a}^i(M)/K_{1})_{\geq n}$  is finite. Since
\begin{center}
$(\Ass_R H_I^i(M)/K_{1})_{\geq n}\cap (\Ass_R H_{\frak a}^i(M))_{\geq n}=\emptyset,$
\end{center}
it follows from Lemma \ref{lem1} that
 \begin{center}
$\cap_{\p \in (\Ass_R H_{\frak a}^i(M))_{\geq n}}\frak p\subsetneqq \cap_{\p \in (\Ass_R H_{\frak a}^i(M)/K_{1})_{\geq n}}\frak p.$
\end{center}

By using the method used in the above, there is a finitely generated submodule $K_{2}/K_{1}$ of $H_{\frak a}^i(M)/K_{1}$ such that

\begin{center}
$(\Ass_R H_{\frak a}^i(M)/K_{2})_{\geq n}\cap (\Ass_R H_{\frak a}^i(M)/K_{1})_{\geq n}=\emptyset,$
\end{center}
and so
\begin{center}
$\cap_{\p \in (\Ass_R H_{\frak a}^i(M)/K_{1})_{\geq n}}\frak p\subsetneqq \cap_{\p \in (\Ass_R H_{\frak a}^i(M)/K_{2})_{\geq n}}\frak p.$
\end{center}
Proceeding in the same way we can find a chain of ideals of $R$,
\begin{center}
$\cap_{\p \in (\Ass_R H_{\frak a}^i(M))_{\geq n}}\frak p\subsetneqq\cap_{\p \in (\Ass_R H_{\frak a}^i(M)/K_{1})_{\geq n}}\frak p\subsetneqq \cap_{\p \in (\Ass_R H_{\frak a}^i(M)/K_{2})_{\geq n}}\frak p\subsetneqq\hdots,$
\end{center}
which is not stable. Consequently, $H_{\frak a}^i(M)$ is in dimension $< n$, as required.\qed \\

Now we are prepared to state and prove the first main result of this paper, which shows that the least integer $i$ such that $H^i_{\frak a}(M)$ is not
in dimension $<n$, equals to $\inf \{f_{\frak aR_{\frak p}}(M_{\frak p})\,\,|\,\,\frak p\in \Supp M/\frak aM\,\,{\rm and} \,\,\dim R/\frak
p\geq n\}.$ This generalizes  the main result of Asadollahi-Naghipour \cite{AN}.

\begin{thm}\label{th1}
Let $R$ be a Noetherian ring, ${\frak a}$ an ideal of $R$ and $M$  a finitely generated $R$-module. Then for any non-negative integer $n$,
$$f_{\frak a}^n(M)={\rm h}_{\frak a}^n(M).$$
\end{thm}

\proof The assertion follows from the definition and Proposition \ref{prop}.\qed\\

As a first application of Theorem 2.10, we  show that the first non-minimax (resp. non-skinny)  local cohomology modules
$H^{i}_{\frak a}(M)$ of a finitely module $M$ over a  Noetherian ring $R$ with respect to an ideal $\frak a$ is
${\rm h}^1_{\frak a}(M)$ (resp. ${\rm h}^2_{\frak a}(M)$).\\
\begin{cor}\label{211}
Let $R$ be a Noetherian ring, ${\frak a}$ an ideal of $R$ and $M$  a finitely generated $R$-module. Then

{\rm (i)}  ${\rm h}^1_{\frak a}(M)=\inf\{i\in \Bbb{N}_0\,\,|\,\,H^i_{\frak a}(M)\,\,{\rm is}\,\,{\rm not}\,\,{\rm minimax}\}.$

{\rm (ii)} ${\rm h}^2_{\frak a}(M)=\inf\{i\in \Bbb{N}_0\,\,|\,\,H^i_{\frak a}(M)\,\,{\rm is}\,\,{\rm not}\,\,{\rm skinny}\},$ whenever $R$
is semilocal.
\end{cor}

\proof $\rm (i)$  follows from \cite[Corollary 2.4]{BNS} and Theorem 2.10.  To prove $\rm (ii)$ use   \cite[Proposition 3.7]{BNS} and Theorem 2.10. \qed\\

\begin{prop}
Let $R$ be a Noetherian ring and $n$ a non-negative integer. Let
 \[0\longrightarrow M^\prime\longrightarrow M\longrightarrow
M^{\prime\prime}\longrightarrow 0\] be an exact sequence of $R$-modules.
Then $M$  is in dimension $< n$ if and only if  $M^\prime$ and
$M^{\prime\prime}$ are both  in dimension $< n$.
\end{prop}
\proof We may suppose for the proof that $M^\prime$ is a submodule
of $M$ and that $M^{\prime\prime}=M/M^\prime$. If $M$ is  in dimension $< n$, then it is easy to verify that
$M^\prime$ and $M/M^\prime$ are  in dimension $< n$. Now, suppose that
$M^\prime$ and $M/M^\prime$ are  in dimension $< n$.
Then there exists a finitely generated
submodule $T$ of $M^{\prime}$ such that $\dim \Supp M^{\prime}/T< n$.
Let $N^{\prime}=M^{\prime}/T$ and $N=M/T$. Then we obtain the exact
sequence $$0 \longrightarrow N^{\prime} \longrightarrow N
\longrightarrow N/N^{\prime} \longrightarrow 0,$$ where $\dim\Supp N^{\prime}< n$ and
$N/N^{\prime}$ is in dimension $< n$, (note that $N/N^{\prime}
\cong M/M^{\prime}$). Now, since $N/N^{\prime}$ is in dimension $< n$ it
follows from the definition that there is a finitely generated
submodule $L/N^{\prime}$ of $N/N^{\prime}$ such that $\dim\Supp N/L< n$. As $L/N^{\prime}$ is finitely generated, it follows that
$L=N^{\prime}+K$ for some finitely generated submodule $K$ of $L$.
Then it follows from $L/K \cong N^{\prime}/K \cap N^{\prime}$ that
$\dim\Supp L/K< n$. Therefore the exact sequence $$0
\longrightarrow L/K \longrightarrow N/K \longrightarrow N/L
\longrightarrow 0$$ implies that $\dim\Supp N/K< n$. Consequently $N$
is  in dimension $< n$. Since $N=M/T$, it follows that $K=S/T$ for some
submodule $S$ of $M$ containing $T$. As $T$ and $K$ are finitely generated, we
deduce that $S$ is also finitely generated. Now Because $$M/S \cong
(M/T)/(S/T) = N/K,$$  it yields that $\dim\Supp M/S <n$, and so by definition, $M$ is
in dimension $< n$ and the claim is true. \qed\\

Before bringing the next results, let us recall that a  full subcategory $\mathcal{S}$ of the category of $R$-modules is
called  a {\it Serre subcategory}, when it is closed under taking submodules, quotients and extensions.
One can easily check that the subcategories of, finitely generated, minimax, skinny, and
Matlis reflexive modules are examples of Serre subcategory. The following result provides a new class of Serre subcategory.
\begin{cor}
For any non-negative integer $n$, the class of in dimension $< n$ modules over a  Noetherian ring $R$ consists  a  Serre subcategory of
the category of $R$-modules.
\end{cor}
\proof The assertion follows immediately from Proposition 2.12. \qed \\
\begin{cor}
 Let $R$ be a Noetherian ring and  $n$ a non-negative integer. Then any quotient
of an  in dimension $< n$ module, as well as any finite direct sum of
 in dimension $<n$ modules, is  in dimension $< n$.
\end{cor}
\proof The assertion follows from  definition and Proposition 2.12.\qed\\

\begin{cor}\label{211}
Let $R$ be a Noetherian ring, ${\frak a}$ an ideal of $R$ and $M$  a
finitely generated $R$-module.  For a non-negative integer $n$, let $t=f_{\frak a}^n(M)$. Then for each submodule $N$ of
$\oplus_{i=0}^{t-1}H_{\frak a}^i(M)$, the set $$(\Ass_R(\oplus_{i=0}^{t-1}H_{\frak a}^i(M)/N))_{\geq n}$$ is finite.
\end{cor}
\proof Since by Proposition \ref{prop},  $\oplus_{i=0}^{t-1}H_{\frak a}^i(M)/N$ is in dimension $< n$,  the
assertion follows from Lemma \ref{lem2}.\qed \\

\begin{cor}
 Let $R$ be a Noetherian ring and $n$ a  non-negative integer.  Let $M, N$ be
$R$-modules such that $M$ is finitely generated  and $N$ is  in dimension $< n$.
Then ${\rm Ext}_R^i(M,N)$ and ${\rm Tor}_i^R(M,N)$ are
in dimension $< n$  modules for all $i$. In particular, for any ideal ${\frak a}$ of $R$, the
$R$-modules ${\rm Ext}_R^i(R/{\frak a}, N)$ and ${\rm Tor}_i^R(R/{\frak a}, N)$ are   in dimension $< n$,  for all $i$.
\end{cor}
\proof As $R$ is Noetherian and $M$ is finitely generated, it
follows that $M$ possesses a free resolution
\[\mathbb{F_\bullet}:\cdots\rightarrow F_s\rightarrow
F_{s-1}\rightarrow\cdots\rightarrow F_1\rightarrow F_0\rightarrow
0,\] whose free modules have finite ranks.

Thus ${\rm Ext}_R^i(M,N)=H^i({\rm Hom}_R(\mathbb{F_\bullet},N))$ is
a subquotient of a direct sum of finitely many copies of $N$.
Therefore, it follows from Corollary 2.14  that ${\rm Ext}_R^i(M,N)$ is  in dimension $< n$  for all $i\geq 0$. By using a
similar proof as above we can deduce that ${\rm Tor}_i^R(M,N)$ is in dimension $< n$  for all $i\geq 0$.\qed\\

The following proposition is needed in the proof of the second main theorem of this section.

\begin{prop}
  Let $R$ be a Noetherian ring and $n$ a  non-negative integer.  Let $M$ be a
finitely generated $R$-module and $N$ an arbitrary $R$-module.
 Let $t$ be a non-negative integer such that ${\rm Ext}_R^i(M,N)$
is in dimension $< n$ for all $i\leq t$. Then for any finitely
  generated $R$-module $L$ with ${\rm Supp}\,L\subseteq {\rm Supp}\,M$, ${\rm Ext}_R^i(L,N)$
   is in dimension $< n$  for all $i\leq t$.
\end{prop}
\proof Since ${\rm Supp}\,L\subseteq {\rm Supp}\,M$, it follows from the
Gruson's Theorem
 (cf. \cite[Theorem 4.1]{Va}), that there exists a chain \[0=L_0\subset L_1\subset\cdots\subset L_k=L,\]
 such that the factors $L_j/L_{j-1}$ are homomorphic images of a direct sum of finitely
 many copies of $M$. Now consider the exact sequences
\[0\rightarrow K\rightarrow
  M^r\rightarrow L_1\rightarrow 0  \ \ \ \ \ \ \ \ \ \ \ \ \ \ \ \ \ \ \      \]

\[0\rightarrow L_1\rightarrow
   L_2\rightarrow L_2/L_1\rightarrow 0  \ \ \ \ \ \ \ \ \ \ \ \ \ \ \ \ \ \ \      \]
                               $$\vdots$$
                                       \[0\rightarrow L_{k-1}\rightarrow
  L_k\rightarrow L_k/L_{k-1}\rightarrow 0,  \ \ \ \ \ \ \ \ \ \ \ \ \ \ \ \ \ \ \  \]

for some positive integer $r$.

Now from the long exact sequence

 \[\cdots\rightarrow {\rm Ext}_R^{i-1}(L_{j-1},N)\rightarrow
   {\rm Ext}_R^{i}(L_j/L_{j-1},N)\rightarrow {\rm
   Ext}_R^{i}(L_j,N) \rightarrow {\rm Ext}_R^i(L_{j-1},N)\rightarrow\cdots,\]

   and an easy induction on $k$, it suffices  to prove the case when
   $k=1$.

 Thus there is an exact sequence \[0\rightarrow K\rightarrow
  M^r\rightarrow L\rightarrow 0  \ \ \ \ \ \ \ \ \ \ \ \ \ \ \ \ \ \ \      (\ast)\] for some $r\in\Bbb N$
and some finitely generated $R$-module $K$.

    Now, we use induction on $t$. First, ${\rm Hom}_R(L,N)$ is a submodule
     of ${\rm Hom}_R(M^r,N)$; hence in view of assumption and
     Corollary 2.14, ${\rm Ext}_R^0(L,N)$ is in dimension $< n$.
     So assume that $t>0$ and that ${\rm Ext}_R^j(L^\prime,N)$
     is in dimension $< n$  for every finitely
  generated $R$-module $L^\prime$ with ${\rm Supp}\,L^\prime\subseteq {\rm Supp}\,M$ and
  all $j\leq t-1$. Now, the exact sequence $(\ast)$ induces the long
  exact sequence \[\cdots\rightarrow {\rm Ext}_R^{i-1}(K,N)\rightarrow
   {\rm Ext}_R^{i}(L,N)\rightarrow {\rm
   Ext}_R^{i}(M^r,N)\rightarrow\cdots,\]so that, by the inductive
   hypothesis, ${\rm Ext}_R^{i-1}(K,N)$ is in dimension $< n$  for all $i\leq
   t$. On the
   other hand, according to Corollary  2.14,
   $$\Ext_R^{i}(M^r,N)\cong\,  \stackrel{r} {\oplus}{\rm Ext}_R^{i}(M,N)$$
   is in dimension $<n$. Therefore, it follows from Proposition 2.12 that ${\rm Ext}_R^{i}(L,N)$ is in dimension $< n$ for all $i\leq t$,
   and this completes  the inductive step.\qed \\

\begin{cor}
 Let $R$ be a Noetherian ring and $\frak a$ an ideal of $R$. Assume that  $t, n$ are non-negative integers. Then, for any $R$-module $M$ the following
conditions are equivalent:

{\rm(i)} ${\rm Ext}_R^i(R/\frak a, M)$ is in dimension $< n$ for all $i\leq t$.

{\rm(ii)} For any ideal $\frak b$ of $R$ with $\frak b\supseteq \frak a$,  ${\rm Ext}_R^i(R/\frak b,M)$ is in dimension $< n$ for all $i\leq t$.

{\rm(iii)} For any finitely generated $R$-module $N$ with ${\rm
Supp}(N)\subseteq V(\frak a)$,  ${\rm Ext}_R^i(N,M)$ is in dimension $< n$ for all $i\leq t$.

{\rm(iv)} For any minimal prime ideal $\frak p$ over $\frak a$, ${\rm Ext}_R^i(R/\frak p,M)$ is in dimension $< n$ for all $i\leq t$.

\end{cor}
\proof In view of the Proposition 2.17, it is enough to show that {\rm(iv)} implies
 {\rm(i)}. To do this,
 let $\frak p_1,\ldots,\frak p_s$ be the minimal primes of $\frak a$. Then, by assumption
  the $R$-modules ${\rm Ext}_R^i(R/\frak p_j,M)$ are in dimension $< n$ for all
   $j=1,2,\dots, s$. Hence by Corollary 2.14, $$\oplus_{j=1}^s{\rm Ext}_R^{i}(R/\frak p_j,M)\cong{\rm
   Ext}_R^{i}(\oplus_{j=1}^sR/\frak p_j,M)$$ is in dimension $< n$. Since
 ${\rm Supp}(\oplus_{j=1}^sR/\frak p_j)={\rm Supp}R/\frak a$, it follows from
Proposition 2.17 that ${\rm Ext}_R^i(R/\frak a, M)$ is in dimension $<n$, as required.\qed\\

 Following we let $\mathcal{S}$ denote a Serre subcategory of the category of $R$-modules.

\begin{lem}
    \label{2.1}
 Let $R$ be a Noetherian ring and $\frak a$ an ideal of $R$. Let $s$  be a
 non-negative integer and $M$  an $R$-module such that $\Ext_R^s(R/\frak a,M)\in \mathcal{S}$.
If $\Ext_R^j(R/\frak a,H^i_\frak a(M))\in \mathcal{S}$ for all $i<s$ and all $j\geq0$, then
$\Hom_R(R/\frak a,H^s_{\frak a}(M))\in \mathcal{S}$.
\end{lem}
\proof See \cite[Theorem 2.2]{AS}.\qed\\

\begin{prop}
    \label{2.2}
Let $R$ be a Noetherian ring and $\frak a$ an ideal of $R$. Let $s$  be a
 non-negative integer and let $M$ be an $R$-module such that $\Ext_R^{s+1}(R/\frak a,M)\in \mathcal{S}$.
If $\Ext_R^j(R/\frak a,H^i_{\frak a}(M))\in \mathcal{S}$ for all $i<s$ and all $j\geq0$, then
$\Ext^1_R(R/{\frak a},H^s_\frak a(M))\in \mathcal{S}$.
\end{prop}
\proof We use induction on $s$. Let $s=0$. Then the exact sequence
$$0\longrightarrow \Gamma_{\frak a}(M)\longrightarrow M\longrightarrow M/\Gamma_{\frak a}(M)\longrightarrow0,$$
induces the exact sequence
$$\Hom_R(R/{\frak a},M/\Gamma_{\frak a}(M))\longrightarrow \Ext^1_R(R/{\frak a}, \Gamma_{\frak a}(M)) \longrightarrow \Ext^1_R(R/{\frak a}, M).$$
 As $\Hom_R(R/{\frak a}, M/\Gamma_{\frak a}(M))=0$ and $\Ext_R^{1}(R/{\frak a},M)\in \mathcal{S}$, it follows that
$\Ext_R^{1}(R/{\frak a},\Gamma_{\frak a}(M))$ is also  in $\mathcal{S}$.

Now suppose inductively that $s>0$ and that the assertion holds for $s-1$. Again using the exact sequence
$$0\longrightarrow \Gamma_{\frak a}(M)\longrightarrow M\longrightarrow M/\Gamma_{\frak a}(M)\longrightarrow0,$$
for all $j\geq0$, we obtain the following exact sequence,
$$\Ext^j_R(R/{\frak a},M)\longrightarrow \Ext^j_R(R/{\frak a}, M/\Gamma_{\frak a}(M)) \longrightarrow \Ext^{j+1}_R(R/{\frak a}, \Gamma_{\frak a}(M)).$$
Now since by assumption, $\Ext^{s+2}_R(R/{\frak a}, \Gamma_{\frak a}(M))$ and $\Ext^{s+1}_R(R/{\frak a},M)$ are in $\mathcal{S}$, it follows that
 $\Ext^{s+1}_R(R/{\frak a}, M/\Gamma_{\frak a}(M))\in\mathcal{S}$. Also, it follows easily  from assumption and \cite[Corollary 2.1.7]{BS}
that $$\Ext_R^j(R/{\frak a},H^i_{\frak a}(M/\Gamma_{\frak a}(M)))\in \mathcal{S}$$ for all $i<s$ and all $j\geq0$. Therefore we may
assume that $\Gamma_{\frak a}(M)=0$.

Next, let $E_R(M)$ denote the injective hull of $M$. Then $\Gamma_{\frak a}(E_R(M))=0$,
and so it follows from the exact sequence
$$0\longrightarrow M\longrightarrow E_R(M)\longrightarrow E_R(M)/M\longrightarrow 0,$$
that $H^{i+1}_{\frak a}(M)\cong H^{i}_{\frak a}(E_R(M)/M)$ for all $i\geq0$. Also, as $\Hom_R(R/{\frak a}, E_R(M))=0$,
it yields that $$\Ext_R^j(R/{\frak a}, E_R(M)/M\cong \Ext_R^{j+1}(R/{\frak a}, M),$$ for all $j\geq0$.
Consequently the $R$-module $E_R(M)/M$ satisfies our condition hypothesis. Thus $$\Ext^1_R(R/{\frak a},H^{s-1}_{\frak a}(E_R(M)/M))\in \mathcal{S},$$
and so the assertion follows from $H^{s}_{\frak a}(M)\cong H^{s-1}_{\frak a}(E_R(M)/M).$ \qed\\

We are now ready to state and prove the second main result of this paper which is a generalization of the main result of Quy \cite[Theorem 3.2]{Qu}  and  Brodmann-Lashgari  \cite[Theorem 2.2]{BL}.

\begin{thm}
Let $R$ be a Noetherian ring, ${\frak a}$ an ideal of $R$ and $M$ a
finitely generated $R$-module. For a non-negative integer $n$, let $t=f_{\frak a}^n(M)$. Then the following statements hold:

$\rm (i)$ The $R$-modules  ${\rm Ext}_R^i(R/{\frak a},H_{\frak a}^i(M))$ are in dimension $<n$ for $i=0,1,\dots, t-1$ and all integers $j$.

$\rm (ii)$ The $R$-modules ${\rm Hom}_R(R/{\frak a},H_{\frak a}^t(M))$ and ${\rm Ext}_R^1(R/{\frak a},H_{\frak a}^t(M))$ are in dimension $< n$.

$\rm (iii)$  For each finitely generated $R$-module $N$ with ${\Supp}(N)\subseteq V({\frak a})$ the $R$-modules ${\rm Ext}_R^j(N,H_{\frak a}^i(M))$ are in dimension
$< n$ for all $i=0,1,\dots, t-1$ and for all integers $j$.

$\rm (iv)$ The set $({\Ass}_R(H_{\frak a}^t(M)))_{\geq n}$ is finite.

$\rm (v)$ The set ${\Ass}_R(\oplus_{i=0}^tH_{\frak a}^i(M)/K)_{\geq n}$ is finite, for any  finitely generated submodule $K$ of
$\oplus_{i=0}^t(H_{\frak a}^i(M)$.
\end{thm}
\proof  $\rm (i)$  This follows immediately from Proposition 2.9 and  Corollary 2.16.  In order to show (ii) use part (i),  Lemma 2.19, Theorem 2.10 and Proposition  2.20.
Part (iii)  follows from  Proposition 2.17 and part (i).  Moreover, (iv) follows readily from Lemma \ref{lem2}, part  (ii) and the fact that
$$\Ass_R\Hom_R(R/\frak a, H_{\frak a}^t(M))= \Ass_RH_{\frak a}^t(M).$$

 Finally, (v) follows easily from  Corollary \ref{211},  part (iv) and the exact sequence
 $$\Hom_{R}(R/\frak a, \oplus_{i=0}^tH_{\frak a}^i(M)) \longrightarrow \Hom_{R}(R/\frak a, \oplus_{i=0}^tH_{\frak a}^i(M)/K) \longrightarrow \Ext^{1}_{R}(R/\frak a, K).$$  \qed \\

The paper ends a  result  about the finiteness Bass numbers of a certain local cohomology modules. Recall that for any prime ideal $\frak p$ of
$R$ and an $R$-module $L$, the $i$-th {\it Bass number}  $\mu^i(\frak p, L)$ is defined to be $\dim_{k(\frak p)}\Ext^i_{R_{\frak p}}(k(\frak p), L_{\frak p})$,
 where $k(\frak p)=R_{\frak p}/\frak pR_{\frak p}$. \\
\begin{cor}
Let $R$ be a Noetherian ring, $\frak a$ an ideal of $R$ and $M$  a finitely generated $R$-module.   For a non-negative integer $n$, let
$t={\rm h}_{\frak a}^n(M)$. Let $\p\in \Spec R$ be such that $\dim R/\frak p\geq n$. Then

$\rm (i)$ The Bass numbers $\mu^j(\p,H_{\frak a}^i(M))$ are finite for all $0\leq i\leq t-1$  and all integers $j$.

$\rm (ii)$ The Bass numbers $\mu^j(\p, H_{\frak a}^t(M))$ are finite for $j=0,1$.
\end{cor}
\proof  The part $\rm (i)$  follows from Theorem 2.10 and the fact that the Bass numbers of finitely generated modules are finite. Also,   $\rm (ii)$   follows from Theorem 2.10  and \cite[Corollary 3.5]{Kh}.\qed \\
%\end{align*}
%\section{The Results}

\begin{center}
{\bf Acknowledgments}
\end{center}
The authors  are deeply grateful to the
referee for his or her valuable suggestions on the paper and for
drawing the authors'  attention to Corollary 2.11. Also, we
 would like to thank Dr.  Kamal Bahmanpour for reading of the first draft and valuable discussions.
Finally, the authors would like to thank Tabriz branch, Islamic Azad University for the financial support of this research, which
is based on a research project contract
% ----------------------------------------------------------------

\end{document}